\documentclass[a4paper]{article}
\usepackage[utf8]{inputenc}
\usepackage{lmodern}
\usepackage[english]{babel}
\usepackage[T1]{fontenc}
\usepackage{amsmath,amsthm,amssymb}
\usepackage[colorlinks=true,allcolors=black]{hyperref}
\usepackage{csquotes}
\usepackage{enumitem}

\usepackage{tikz}
\usepackage{pgfplots}
\pgfplotsset{compat=1.18}
\usetikzlibrary{datavisualization}
\usetikzlibrary{datavisualization.formats.functions}
\usetikzlibrary{matrix, calc, decorations}

\usepackage[backend=biber]{biblatex}

\usepackage[left=1in,right=1in,top=1in,bottom=.5in,includeheadfoot]{geometry}
\usepackage[colorlinks=true,allcolors=black]{hyperref}

\title{Linear Equations for Solving Partial~Latin~Squares and Sudokus}
\author{Ralf Pöppel \\ \url{mailto:ralf@poeppel-familie.de}}
\date{2025-08-25}
\begin{document}
\maketitle

\begin{abstract}
For Latin squares we know the units: rows and columns, have fixed sum. The same holds for rows and columns
and the units blocks of Sudokus.

The summation of the elements of a unit is a linear equation. So the set of equations of all units
gives a system of linear equations that models the relations of the variables of
Latin squares and Sudoku puzzles.

Every completed Latin square or Sudoku satisfies this system, enabling the analysis and completion of partially filled puzzles.

The linear system provides a necessary condition for a solution, but not a sufficient one.
Thus solving the system does not guarantee the completion in all cases. In this article we analyze
this linear system and its limitation in solving puzzles.

We show that for Latin squares of size $n \times n$ with $k$ unknowns, and for Sudokus of size $n \times n$ with $n = m \cdot l$ and $k$ unknowns, the system has full rank only if:
\[
k \leq 2n - 1 \quad \text{(Latin squares)}, \quad k \leq 2n - 1 + (l - 1)(m - 1) \quad \text{(Sudokus)}.
\]
These rank conditions are necessary for the system to uniquely determine the unknown entries under the
linear model. On the other hand, we give an example of a Sudoku with 2 solutions, where this condition holds.

Addionally the linear system captures only one of the three relations of Latin squares and Sudokus, known to
puzzle solvers. To fully solve any puzzle, it must be extended with additional not necessarily linear equations
reflecting the additonal relations.

This article examines the conditions under which Latin squares and Sudoku puzzles can be described and solved as a linear
system of equations.

\end{abstract}

\section{Introduction}

\subsection{Motivation}

In recent years, Sudokus \cite{SudokuSerious,Delahaye2006,wikipedia_en_2025_Sudoku}
have gained increasing popularity as mathematical puzzles.
The task of a Sudoku is to complete a partially filled square grid of $9 \times 9$ numbers in such a way that each natural number
from $1$ to $9$ appears exactly once
in each row, each column, and each $3 \times 3$ subgrid, called a block.
The procedures used when solving Sudokus resemble, to some extent, the back-substitution method used in solving linear systems of equations.
This article examines the conditions under which a Sudoku formulated as a
system of linear equations has a unique solution.

\subsection{Related Work}

There are many articles and books on Sudoku. Most publications present different approaches to
solving Sudokus. Rosenhouse and Taalman \cite{SudokuSerious} provide an overview of various approaches. All differ from the approach presented here.
Chen and Cooper \cite{Chen2017} state:
\begin{quote}
There are numerous books and articles on Sudoku. We do not find any
deterministic algorithms to solve all puzzles. We shall study the structure
and properties of Sudoku puzzles and establish some strategies for solving
puzzles deterministically, i.e. without trial-and-error.
\end{quote}
They provide several logical strategies for solving Sudokus.
A similar approach is used by
Berthier \cite{berthier2007}, who provides 22 logical rules in detail for solving Sudokus and uses an inference machine for solutions.

Arnold et al.\ \cite{Arnold_Sudoku_2010} provide a
system of linear and nonlinear equations over $\mathbb{N}$ and attempted to find a Gröbner basis, with limited success. The employed library algorithm was
able to compute a Gröbner basis only for Sudokus of order 4. Nacin \cite{Nacin2024} presents a system of linear equations for cyclic groups
$\mathbb{Z}_n$ and uses the equations to solve the puzzle step by step.

\subsection{Used Approach}

A simpler task that can be solved analogously to a Sudoku is the completion of a partially filled Latin square \cite{EMS_LatinSquare,Mathworld2025:LatinSquare}.
Latin squares, like Sudokus, are square grids in which the numbers $1$ through $n$ may occur
only once in each row and each column. Every Sudoku is also a Latin square \cite{berthier2007,Nacin2024}.
In other words: a Sudoku is a specialization of a Latin square.
As far as we know, this is the first time Latin squares are analyzed together with Sudokus in this way.

In this text, we provide a linear system of equations over $\mathbb{N}$ derived from the known properties of Latin squares and Sudokus,
and show that the system of equations of a Sudoku is an extension of the system of Latin squares.
We solve the system of equations holistically
using row reduction, also called Gaussian elimination \cite{wikipedia_en_2025_gaussian-elimination,Strang2003}.
This will prove the existence
of partial Latin squares and Sudokus, which can be solved in polynomial time. Gaussian elimination has arithmetic complexity
$\mathcal{O}(n^3)$ \cite{Meier2020}.
Furthermore, we analyze the
feasibility and limitations of this approach for solving Sudokus and completing Latin squares.

We will first examine the conditions under which partially filled Latin squares can be
formulated and solved as a linear system of equations, followed by a consideration of Sudokus. This order
was chosen because Latin squares are generalizations of Sudokus and therefore fewer conditions need to be
considered when formulating Latin squares. This results in smaller systems of equations.

\subsection{Latin Square}
A Latin square \cite{EMS_LatinSquare, Nacin2024} is a square grid with $n$ rows and $n$ columns.
The Latin square is filled with $n$ symbols in such a way that each symbol occurs only once in each row and column.
The number $n$ is called the order of the Latin square.
The set of symbols may consist of numbers, letters, colors, or any other collection of distinct elements.
In many cases, including in this article, the natural numbers from $1$ to $n$ are used as symbols.

The task of completing a partially filled Latin square coresponds to the task of solving a Sudoku.
\autoref{lq4full} shows a partially filled Latin square (the puzzle) and the corresponding fully filled Latin square (the solution) of order 4.

\begin{figure}
\begin{tikzpicture}[scale=0.8]
\node at (-5, 0) { };
\draw[help lines] (0,0) grid (4,4);
    \draw[line width=1.5pt]  (0,0) -- (4,0) -- (4,4) --(0,4) -- (0,0) -- (4,0);
\node  at (0.5, 3.5) {1};
\node  at (2.5, 2.5) {1};
\node  at (3.5, 1.5) {1};
\node  at (3.5, 3.5) {2};
\node  at (0.5, 2.5) {3};
\node  at (1.5, 1.5) {3};
\node  at (3.5, 0.5) {3};
\color{black}
\node  at (1.5, 3.5) {4};
\node  at (2.5, 0.5) {4};
\draw[help lines] (6,0) grid (10,4);
\draw[line width=1.5pt]  (6,0) -- (10,0) -- (10,4) -- (6,4) -- (6,0) -- (10,0);
\node  at (6.5, 3.5) {1};
\node  at (8.5, 2.5) {1};
\node  at (9.5, 1.5) {1};
\node  at (7.5, 0.5) {1};
\node  at (9.5, 3.5) {2};
\node  at (7.5, 2.5) {2};
\node  at (8.5, 1.5) {2};
\node  at (6.5, 0.5) {2};
\node  at (8.5, 3.5) {3};
\node  at (6.5, 2.5) {3};
\node  at (7.5, 1.5) {3};
\node  at (9.5, 0.5) {3};
\color{black}
\node  at (7.5, 3.5) {4};
\node  at (9.5, 2.5) {4};
\node  at (6.5, 1.5) {4};
\node  at (8.5, 0.5) {4};
\end{tikzpicture}
\caption{Example of a Latin square of order $n=4$ \label{lq4full}}
\end{figure}

\subsection{Sudoku}
A Sudoku, like a Latin square, is a square grid with $n$ rows and $n$ columns.
Just like a Latin square, a Sudoku is filled with $n$ symbols such that each symbol appears exactly once in every row and every column.
Usually, Sudokus use only the natural numbers $1, \ldots , n$.

In addition to Latin squares, Sudokus are divided into blocks (sub-squares or sub-rectangles)
of size $l \times m = n$, and the additional condition holds
that each symbol must appear exactly once in every block. Here, $l$ is the number of rows and $m$ the number of columns in a block.

A subset of Latin squares whose order $n$ is not a prime number also satisfies
this additional Sudoku condition.
The most common Sudokus have order $n = 9$ with blocks of size $l \times m = 3 \times 3$.

\autoref{sdk4full} shows an example of a Sudoku of order $n = 4$ with blocks of size $2 \times 2$ as
a puzzle and its solution.

\begin{figure}[h]
\begin{tikzpicture}[scale=0.8]
\node at (-5, 0) { };
\draw[help lines] (0,0) grid (4,4);
\draw[line width=1.5pt]  (0,0) -- (4,0) -- (4,4) -- (0,4) -- (0,0) -- (4,0);
\draw[line width=1.5pt]  (0,2) -- (4,2);
\draw[line width=1.5pt]  (2,0) -- (2,4);
\node  at (0.5, 3.5) {1};
\node  at (2.5, 2.5) {1};
\node  at (3.5, 1.5) {1};
\node  at (3.5, 3.5) {2};
\node  at (0.5, 2.5) {3};
\node  at (1.5, 1.5) {3};
\node  at (3.5, 0.5) {3};
\color{black}
\node  at (1.5, 3.5) {4};
\node  at (2.5, 0.5) {4};
\draw[help lines] (6,0) grid (10,4);
\draw[line width=1.5pt]  (6,0) -- (10,0) -- (10,4) -- (6,4) -- (6,0) -- (10,0);
\draw[line width=1.5pt]  (6,2) -- (10,2);
\draw[line width=1.5pt]  (8,0) -- (8,4);
\node  at (6.5, 3.5) {1};
\node  at (8.5, 2.5) {1};
\node  at (9.5, 1.5) {1};
\node  at (7.5, 0.5) {1};
\node  at (9.5, 3.5) {2};
\node  at (7.5, 2.5) {2};
\node  at (8.5, 1.5) {2};
\node  at (6.5, 0.5) {2};
\node  at (8.5, 3.5) {3};
\node  at (6.5, 2.5) {3};
\node  at (7.5, 1.5) {3};
\node  at (9.5, 0.5) {3};
\color{black}
\node  at (7.5, 3.5) {4};
\node  at (9.5, 2.5) {4};
\node  at (6.5, 1.5) {4};
\node  at (8.5, 0.5) {4};
\end{tikzpicture}
\caption{Example of a Sudoku of order $n = 4, l = 2, m = 2$ \label{sdk4full}}
\end{figure}

Please note the puzzle in \autoref{sdk4full} equals the puzzle in \autoref{lq4full}. This Sudoku can be solved using the rules of Latin Squares only.

\section{Linear System of Equations}
The first step in setting up a linear system of equations \cite{Strang2003} is defining the variables.

The variables for the numbers in the cells of the square grid are denoted by $a_{i,j}$, where $i, j = 1, \ldots, n$.
We have $a_{i,j} \in \{x \in \mathbb{N} \mid 1 \leq x \leq n\}$, where $i$ is the row index and $j$ the column index.
The variables can be written as a square matrix $A$:
\begin{equation*}
A := (a_{i,j}) =
\begin{pmatrix}
a_{1,1} & \ldots & a_{1, n}\\
\vdots & & \vdots \\
a_{n,1} & \ldots &a_{n,n}
\end{pmatrix} .
\end{equation*}
If the variables $a_{i,j}$  are arranged line by line, the result is the vector
\begin{equation}
x = (x_k) = (a_{1,1}, \ldots , a_{1, n},  \quad  \ldots, \quad    a_{n, 1}, \ldots , a_{n, n}) \label{vec_x}, \quad k = 1 \ldots n^2
\end{equation}
with the bijective linear mapping $\varphi$
\begin{eqnarray*}
\varphi : \mathbb{N}^{n \times n} & \rightarrow & \mathbb{N}^{n^2}\\
(a_{i,j}) & \mapsto & (x_k), \text{ mit } \\
k  & = & (i - 1) \cdot n + j,\\
\end{eqnarray*}
and the inverse function
\begin{eqnarray*}
\varphi^{-1} : \mathbb{N}^{n^2} & \rightarrow & \mathbb{N}^{n \times n}\\
(x_k) & \mapsto & (a_{i,j}), \text{ mit }  \\
j & = & 1 + (k - 1) \bmod{n}\\
i & = & (k - j) / n + 1.\\
\end{eqnarray*}

Example of $\varphi$ for $n = 2$:
\begin{equation*}
\begin{array}{c | c|c|c|c}
& a_{1,1} & a_{1,2} & a_{2,1} & a_{2,2}\\
\hline
& x_1 & x_2 & x_3 & x_4\\
\hline
i & 1 & 1 & 2 & 2\\
\hline
j & 1 & 2 & 1 & 2\\
\hline
k & 1 & 2 & 3 & 4\\
\end{array}
\end{equation*}.

\subsection{Latin Square}
\subsubsection{Equations}

In a Latin square with $a_{i,j} \in \{x \in \mathbb{N} \mid 1\leq x \leq n\}$ each unit, either row or column, has the sum
\begin{equation*}
s := \sum_{k=1}^{n} k = \frac{n(n+1)}{2}.
\end{equation*}

This means that the following linear system of equations with $2n$ equations for the $2n$ units,
the $n$ rows and $n$ columns, applies to a Latin square of
order $n$:
\begin{eqnarray}
\sum_{i=1}^{n} a_{i,j} = s, &  j = 1, \ldots , n, & \mbox{ equations $1$ to $n$ for $n$ columns}, \label{eq:lq:a}\\
\sum_{j=1}^{n} a_{i,j} = s,  & i = 1, \ldots , n,  & \mbox{ equations $n + 1$ to $2n$ for $n$ rows}, \nonumber
\end{eqnarray}

and the following system of equations for the variables $x_i$ of the vector $x$:
\begin{eqnarray}
\sum_{i=0}^{n-1} x_{n \cdot i + k} = s, &  k = 1 , \ldots , n, & \mbox{equations $1 , \ldots , n$ for $n$ columns}, \label{eq:lq:x}\\
\sum_{j=1}^{n} x_{n \cdot k + j} = s, & k = 0 , \ldots , n-1,  & \mbox{equations $n + 1 , \ldots , 2n$ for $n$ rows}. \nonumber
\end{eqnarray}

\subsubsection{Transformation to Row Reduced Echelon form}
The system of equations (\ref{eq:lq:x}) is written as a matrix $B x = c$ with the coefficient matrix $B = (b_{i,j})$ and the
vector $c = (c_i)$ on the right-hand side. The coefficient matrix has $2n$ rows and $n^2$ columns and only
entries from the set $\{0, 1\}$. The extended coefficient matrix $(B | c) = (b_{i,j} | c_i)$ is as follows:

\begin{equation*}
\begin{array}{llcl}
i = 1, \ldots , n : & & &\\
&  b_{i,j} =  \left \{ \begin{array}{ll}
                     		1,    & j = i + (k - 1) \cdot n, k = 1 , \ldots , n \\  % k = 1 2 3 , j = i, n+i, (k-1)n + i
                     		0,    & \mbox{otherwise}
                                    \end{array}
                       \right\}  & | & c_i = s ,\\
i = n+1, \ldots , 2n : & & & \\
& b_{i,j} = \left \{ \begin{array}{ll}
                     		1,    & j = (i - 1 - n) \cdot n + k, k = 1 , \ldots , n \\ % k = 1 2 3, j = 1 2 3 + (i-1)n
                     		0,    & \mbox{otherwise}
                                    \end{array}
                     \right\} & | & c_i = s .
\end{array}
\end{equation*}

The lines $i = n+1$ to $i = 2n$ can be split into 2 sets as follows:
\begin{equation*}
\begin{array}{llcl}
i = n+1: & & & \\
& b_{i,j} = \left \{ \begin{array}{ll}
                     		1    &  j = 1 , \ldots , n \\
                     		0    & \mbox{otherwise}
                                    \end{array}
                     \right\} & | & s ,\\
i = n+2, \ldots , 2n : & & & \\
& b_{i,j} = \left \{ \begin{array}{ll}
                     		1    & j = (i - 1 - n) \cdot n + k, k = 1 , \ldots , n \\ % k = 1 2 3, j = 1 2 3 + (i-1)n
                     		0    & \mbox{otherwise}
                                    \end{array}
                     \right\} & | & s.\\
\end{array}
\end{equation*}

If lines $1$ to $n$ remain unchanged and the last n lines of $(B|c)$ are reordered so that the
$n + 1$ line becomes the last line and the lines $n + 2$ to $2n$ retain their order, we get
this system of equations with the lines renumbered:
\begin{equation*}
\begin{array}{llcl}
i = 1, \ldots , n : & & &\\
& b_{i,j} = \left \{ \begin{array}{ll}
                     		1,    & j = i + (k - 1) \cdot n, k = 1 , \ldots , n \\  % k = 1 2 3 , j = i, n+i, (k-1)n + i
                     		0,    & \mbox{otherwise}
                                    \end{array}
                     \right\} & | & s , \\
i = n+1, \ldots , 2n - 1 : & & & \\
& b_{i,j} = \left \{ \begin{array}{ll}
                     		1    & j = (i - n) \cdot n + k, k = 1 , \ldots , n \\ % k = 1 2 3, j = 1 2 3 + (i-1)n
                     		0    & \mbox{otherwise}
                                    \end{array}
                     \right\} & | & s ,\\
i = 2n : & & & \\
& b_{i,j} = \left \{ \begin{array}{ll}
                     		1    &  j = 1 , \ldots , n \\
                     		0    & \mbox{otherwise}
                                    \end{array}
                     \right\} & | & s ,\\
\end{array}
\end{equation*}

The resulting matrix can be transformed to row reduced echolon form in 2 steps.
In the first step, the lines $i = 1, \ldots ,n$ are subtracted from the last line $i = 2n$.

\begin{equation*}
\begin{array}{llcl}
i = 1, \ldots , n : & & &\\
& b_{i,j} = \left \{ \begin{array}{ll}
                     		1,    & j = i + (k - 1) \cdot n, k = 1 , \ldots , n \\  % k = 1 2 3 , j = i, n+i, (k-1)n + i
                     		0,    & \mbox{otherwise}
                                    \end{array}
                     \right\} & | & s , \\
i = n+1, \ldots , 2n - 1 : & & & \\
& b_{i,j} = \left \{ \begin{array}{ll}
                     		1,    & j = (i - n) \cdot n + k, k = 1 , \ldots , n \\ % k = 1 2 3, j = 1 2 3 + (i-1)n
                     		0,    & \mbox{otherwise}
                                    \end{array}
                     \right\} & | & s, \\
i = 2n : & & & \\
& b_{i,j} = \left \{ \begin{array}{ll}
                     		-1   & j = n + 1, \ldots, n^2 \\
                     		0    & \mbox{otherwise}
                                    \end{array}
                     \right\} & | & - s\cdot (n-1) .\\
\end{array}
\end{equation*}

In the second step, the lines $i = (n+1), \ldots , (2n -1)$ are added to the last line.
\begin{equation*}
\begin{array}{llcl}
i = 1, \ldots , n : & & &\\
& b_{i,j} = \left \{ \begin{array}{ll}
                     		1    & j = i + (k - 1) \cdot n, k = 1 , \ldots , n \\  % k = 1 2 3 , j = i, n+i, (k-1)n + i
                     		0    & \mbox{otherwise}
                                    \end{array}
                     \right\} & | & s, \\
i = n+1, \ldots , 2n - 1 : & & & \\
& b_{i,j} = \left \{ \begin{array}{ll}
                     		1    & j = (i - n) \cdot n + k, k = 1 , \ldots , n \\ % k = 1 2 3, j = 1 2 3 + (i-1)n
                     		0    & \mbox{otherwise}
                                    \end{array}
                     \right\} & | & s, \\
i = 2n : & & & \\
& b_{i,j} = \left \{ \begin{array}{ll}
				0,  & j = 1 , \ldots , n^2
                                   \end{array}
                     \right\} & | & 0 .
\end{array}
\end{equation*}

With this transformation, the last row becomes a zero row and the matrix is in row reduced echolon form. Obviously,
the $n+1$st row of the system of equations (\ref{eq:lq:x}) is linearly dependent on the other rows.
Now let us consider the linear combination of the rows. Let the rows of $(B | c)$ be $r_i, i = 1 \ldots 2n$. Then we get:
\begin{equation*}
r_{n+1} - \sum_{i=1}^n r_i + \sum_{i = n+2}^{2n} r_i
\end{equation*}
Inserting the equations of the rows gives:
\begin{equation*}
\sum_{j=1}^{n} x_{n \cdot j}  -  \sum_{k = 1}^n \sum_{i=0}^{n-1} x_{n \cdot i + k} + \sum_{k=1}^{n-1}\sum_{j=1}^{n} x_{n \cdot k + j}
 =  s - n \cdot s + (n -1) \cdot s.
\end{equation*}
After renaming the indices
\begin{equation*}
\sum_{j=1}^{n} x_{n \cdot j}  -  \sum_{j = 1}^n \sum_{i=0}^{n-1} x_{n \cdot i + j} + \sum_{i=1}^{n-1}\sum_{j=1}^{n} x_{n \cdot i + j} =  0,
\end{equation*}
Swapping the sums
\begin{equation*}
\sum_{j=1}^{n} x_{n \cdot j}  -  \sum_{i=0}^{n-1} \sum_{j = 1}^n x_{n \cdot i + j} + \sum_{i=1}^{n-1}\sum_{j=1}^{n} x_{n \cdot i + j} =  0
\end{equation*}
and splitting the 2nd sum follows:
\begin{equation*}
\sum_{j=1}^{n} x_{n \cdot j}  - \sum_{j=1}^{n} x_{n \cdot j} -  \sum_{i=1}^{n-1} \sum_{j = 1}^n x_{n \cdot i + j} + \sum_{i=1}^{n-1}\sum_{j=1}^{n} x_{n \cdot i + j} =  0.
\end{equation*}
This reduces the system of equations to the following $2n - 1$ equations in row reduced echolon form:
\begin{eqnarray*}
\sum_{i=0}^{n-1} x_{n \cdot i + k} = s, &  k = 1 , \ldots , n, & \mbox{ equations 1 to $n$ for $n$ columns}, \\
\sum_{j=1}^{n} x_{n \cdot k + j} = s, & k = 1 , \ldots , n-1,  & \mbox{ equations $n + 1$ to $2n - 1$ for $n - 1$ rows}. \\
\end{eqnarray*}
The equations for the variables $a_{i,j}$ are as follows:
\begin{eqnarray}
\sum_{i=1}^{n} a_{i,j} = s, &  j = 1, \ldots , n, & \mbox{ equations 1 to $n$ for $n$ columns}, \label{eq:zeilenstufen} \\
\sum_{j=1}^{n} a_{i,j} = s,  & i = 1, \ldots , n-1,  & \mbox{ equations $n + 1$ to $2n - 1$ for $n- 1$ rows}.\nonumber
\end{eqnarray}

As the system of equations is in row reduced echolon form, the remaining equations are linearly independent. This system
of equations has the rank $2n - 1$. This means that puzzles of Latin squares of order~$n$ can be solved with this equations having a
maximum of $2n - 1$ unknowns.

Each line of the system of equations can only be used to determine one unknown. Each equation contains
exactly $n$ entries not equal to $0$, which means that a maximum of $n - 1$ variables per equation must be
given. If more variables are given per equation, one equation is superfluous and the rank is reduced by~1.

If all variables are unknown, there are $n!$ possibilities for a valid assignment of the variables. This
means that the number of possible Latin squares increases at most with the factorial of the order and thus
disproportionately.

If the pivot variables of the rows of the coefficient matrix $(B | c)$ in row reduced echolon form are selected as unknowns
(underlined in Fig. 3), this results in a Latin square in which the unknowns are in the 1st row and the 1st
column.

\begin{figure}[ht]
\begin{tikzpicture}[scale=0.8]
\node at (-1, 0) { };
\draw[help lines] (6,0) grid (10,4);
\draw[line width=1.5pt]  (6,0) -- (10,0) -- (10,4) -- (6,4) -- (6,0) -- (10,0);
\node  at (6.5, 3.5) {\underline{$a_{1,1}$}};
\node  at (7.5, 3.5) {\underline{$a_{1,2}$}};
\node  at (8.5, 3.5) {\ldots};
\node  at (9.5, 3.5) {\underline{$a_{1,n}$}};
\node  at (6.5, 2.5) {\underline{$a_{2,1}$}};
\node  at (6.5, 1.5) {\vdots};
\node  at (6.5, 0.5) {\underline{$a_{n,1}$}};
\color{black}
\node  at (8.5, 2.5) {\ldots};
\node  at (9.5, 1.5) {\vdots};
\node  at (7.5, 0.5) {$a_{n,2}$};
\node  at (7.5, 2.5) {$a_{2,2}$};
\node  at (7.5, 1.5) {\vdots};
\node  at (9.5, 0.5) {$a_{n,n}$};
\node  at (9.5, 2.5) {$a_{2,n}$};
\node  at (8.5, 0.5) {\ldots};
\end{tikzpicture}
\caption{Partial Latin square of order $n$  \label {lqnpart}}
\end{figure}
This representation also shows that a Latin square of order $n \geq 3$ with 4 unknowns in 2 rows and 2
columns can have one or two solutions, while a Latin square with 3 unknowns of order $n \geq 2$ always
has a unique solution.

For each latin square, we can create this linear system of equations and transform it into row reduced
echolon form and compute the rank. If the rank equals the number of unknows we have a unique solution.
In this case a second check is required. Only if all values of the variables are in $\{1, \ldots, n\}$ we
have a valid Latin square.
If the rank is smaller the puzzle has either more as one solution or exhibits a contradiction.

Overall, only Latin squares with 3 up to a maximum of $2n - 1$ unknowns can therefore be represented and solved as a
linear system of equations.

\qed

\subsubsection{Example}

To illustrate these relationships, let us consider the Latin square in \autoref{lq4full}.
Having $a_{i,j} \in \{1, 2 , 3, 4\}$ and $s = 10$ the system of equations with 8 equations reads:
\begin{eqnarray*}
a_{1,1} + a_{2,1} + a_{3, 1} + a_{4,1} & = & 10\\
a_{1,2} + a_{2,2} + a_{3, 2} + a_{4,2} & = & 10\\
a_{1,3} + a_{2,3} + a_{3, 3} + a_{4,3} & = & 10\\
a_{1,4} + a_{2,4} + a_{3, 4} + a_{4,4} & = & 10\\
\\
a_{1,1} + a_{1,2} + a_{1, 3} + a_{1,4} & = & 10\\
a_{2,1} + a_{2,2} + a_{2, 3} + a_{2,4} & = & 10\\
a_{3,1} + a_{3,2} + a_{3, 3} + a_{3,4} & = & 10\\
a_{4,1} + a_{4,2} + a_{4, 3} + a_{4,4} & = & 10 . \\
\end{eqnarray*}
The vector x results in the Equations:
\begin{equation}
\begin{array}{rcl}
x_1 + x_5 + x_9 + x_{13} & = & 10\\
x_2 + x_6 + x_{10} + x_{14} & = & 10\\
x_3 + x_7 + x_{11} + x_{15} & = & 10\\
x_4 + x_8 + x_{12} + x_{16} & = & 10\\
\\
x_1 + x_2 + x_3 + x_4 & = & 10\\
x_5 + x_6 + x_7 + x_8 & = &  10\\
x_9 + x_{10} + x_{11} + x_{12} &  = & 10\\
x_{13} + x_{14} + x_{15} + x_{16} & = & 10 . \\
\end{array}
\label{lq4xall}
\end{equation}

The system of equations can now be specified in matrix notation with the extended coefficient matrix~$(B|c)$,
where the 5th row is moved to the end,
\begin{equation}
(B | c) =
\left(
\begin{array}{rrrrcrrrrcrrrrcrrrr}
1 & 0 & 0 & 0 & & 1 & 0 & 0 & 0 & & 1 & 0 & 0 & 0 && 1 & 0 & 0 & 0\\
0 & 1 & 0 & 0 & & 0 & 1 & 0 & 0 & & 0 & 1 & 0 & 0 && 0 & 1 & 0 & 0\\
0 & 0 & 1 & 0 & & 0 & 0 & 1 & 0 & & 0 & 0 & 1 & 0 && 0 & 0 & 1 & 0\\
0 & 0 & 0 & 1 & & 0 & 0 & 0 & 1 & & 0 & 0 & 0 & 1 && 0 & 0 & 0 & 1\\
\\
0 & 0 & 0 & 0 & & 1 & 1 & 1 & 1 & & 0 & 0 & 0 & 0 && 0 & 0 & 0 & 0\\
0 & 0 & 0 & 0 & & 0 & 0 & 0 & 0 & & 1 & 1 & 1 & 1 && 0 & 0 & 0 & 0\\
0 & 0 & 0 & 0 & & 0 & 0 & 0 & 0 & & 0 & 0 & 0 & 0 && 1 & 1 & 1 & 1\\
\\
1 & 1 & 1 & 1 &  & 0 & 0 & 0 & 0 & & 0 & 0 & 0 & 0 && 0 & 0 & 0 & 0\\
\end{array}
\right.
\left|
\begin{array}{r}
10\\
10\\
10\\
10\\
\\
10\\
10\\
10\\
\\
10\\
\end{array}
\right) .
\end{equation}

If you then subtract the first 4 rows from the 8th row and add rows 5 to 7 to the 8th row, the matrix is
transformed to row reduced echolon form

\begin{equation}
(B|c) =
\left(
\begin{array}{rrrrcrrrrcrrrrcrrrr}
1 & 0 & 0 & 0 && 1 & 0 & 0 & 0 && 1 & 0 & 0 & 0 && 1 & 0 & 0 & 0\\
0 & 1 & 0 & 0 && 0 & 1 & 0 & 0 && 0 & 1 & 0 & 0 && 0 & 1 & 0 & 0\\
0 & 0 & 1 & 0 && 0 & 0 & 1 & 0 && 0 & 0 & 1 & 0 && 0 & 0 & 1 & 0\\
0 & 0 & 0 & 1 && 0 & 0 & 0 & 1 && 0 & 0 & 0 & 1 && 0 & 0 & 0 & 1\\
\\
0 & 0 & 0 & 0 && 1 & 1 & 1 & 1 && 0 & 0 & 0 & 0 && 0 & 0 & 0 & 0\\
0 & 0 & 0 & 0 && 0 & 0 & 0 & 0 && 1 & 1 & 1 & 1 && 0 & 0 & 0 & 0\\
0 & 0 & 0 & 0 && 0 & 0 & 0 & 0 && 0 & 0 & 0 & 0 && 1 & 1 & 1 & 1\\
\\
0 & 0 & 0 & 0 && 0 & 0 & 0 & 0 && 0 & 0 & 0 & 0 && 0 & 0 & 0 & 0\\
\end{array}
\right.
\left|
\begin{array}{r}
10\\
10\\
10\\
10\\
\\
10\\
10\\
10\\
\\
0\\
\end{array}
\right),
\label{lq4koeff}
\end{equation}
and the 8th line becomes the zero line.

The rank of the matrix can be read off directly here as 7. This means that Latin squares of order 4 can be
completed with a maximum of 7 unknowns. Each unit can be used to determine an unknown if a maximum
of 3 variables per unit of the coefficient matrix are known. If 4 variables per unit are known, one equation
is omitted and the rank of the matrix B is reduced by 1.
When choosing the unknowns, it should be noted that each number may only occur once per row and
per column.

A partial Latin square of order 4 with rank 7 is created if, starting from the complete Latin square (Fig. ~\ref{lq4full}),
the pivot elements of the non-zero rows of the linear equation system are taken as dependent variables and the
remaining variables are set as known. This procedure turns a complete Latin square into a partial Latin
square with rank 7 (see Fig. \ref{lq4rg7}).
\begin{figure}[h]
\begin{tikzpicture}[scale=0.8]
\node at (-5, 0) { };
\draw[help lines] (0,0) grid (4,4);
    \draw[line width=1.5pt]  (0,0) -- (4,0) -- (4,4) -- (0,4) -- (0,0) -- (4,0);
\node  at (0.5, 3.5) {$\not 1$};
\node  at (2.5, 2.5) {1};
\node  at (3.5, 1.5) {1};
\node  at (1.5, 0.5) {1};
\node  at (3.5, 3.5) {$\not 2$};
\node  at (1.5, 2.5) {2};
\node  at (2.5, 1.5) {2};
\node  at (0.5, 0.5) {$\not 2$};
\node  at (2.5, 3.5) {$\not 3$};
\node  at (0.5, 2.5) {$\not 3$};
\node  at (1.5, 1.5) {3};
\node  at (3.5, 0.5) {3};
\color{black}
\node  at (1.5, 3.5) {$\not 4$};
\node  at (3.5, 2.5) {4};
\node  at (0.5, 1.5) {$\not 4$};
\node  at (2.5, 0.5) {4};
\draw[help lines] (6,0) grid (10,4);
    \draw[line width=1.5pt]  (6,0) -- (10,0) -- (10,4) -- (6,4) -- (6,0) -- (10,0);
\node  at (8.5, 2.5) {1};
\node  at (9.5, 1.5) {1};
\node  at (7.5, 0.5) {1};
\node  at (7.5, 2.5) {2};
\node  at (8.5, 1.5) {2};
\node  at (7.5, 1.5) {3};
\node  at (9.5, 0.5) {3};
\color{black}
\node  at (9.5, 2.5) {4};
\node  at (8.5, 0.5) {4};
\end{tikzpicture}
\caption{Generation of partial Latin square of order n= 4 with coefficient matrix of rank 7 \label{lq4rg7}}
\end{figure}

Adding additional equations for the known variables to the system of equations (\ref{lq4xall}),
rearranging the rows
accordingly and computing the row reduced echolon form, results in
the following extended coefficient matrix with full rank:
\begin{equation}
\left(
\begin{array}{rrrrcrrrrcrrrrcrrrr}
1 & 0 & 0 & 0 && 0 & -1 & -1 & -1 && 0 & -1 & -1 & -1 && 0 & -1 & -1 & -1\\
0 & 1 & 0 & 0 && 0 & 1 & 0 & 0 && 0 & 1 & 0 & 0 && 0 & 1 & 0 & 0\\
0 & 0 & 1 & 0 && 0 & 0 & 1 & 0 && 0 & 0 & 1 & 0 && 0 & 0 & 1 & 0\\
0 & 0 & 0 & 1 && 0 & 0 & 0 & 1 && 0 & 0 & 0 & 1 && 0 & 0 & 0 & 1\\
\\
0 & 0 & 0 & 0 && 1 & 1 & 1 & 1 && 0 & 0 & 0 & 0 && 0 & 0 & 0 & 0\\
0 & 0 & 0 & 0 && 0 & 1 & 0 & 0 && 0 & 0 & 0 & 0 && 0 & 0 & 0 & 0\\
0 & 0 & 0 & 0 && 0 & 0 & 1 & 0 && 0 & 0 & 0 & 0 && 0 & 0 & 0 & 0\\
0 & 0 & 0 & 0 && 0 & 0 & 0 & 1 && 0 & 0 & 0 & 0 && 0 & 0 & 0 & 0\\
\\
0 & 0 & 0 & 0 && 0 & 0 & 0 & 0 && 1 & 1 & 1 & 1 && 0 & 0 & 0 & 0\\
0 & 0 & 0 & 0 && 0 & 0 & 0 & 0 && 0 & 1 & 0 & 0 && 0 & 0 & 0 & 0\\
0 & 0 & 0 & 0 && 0 & 0 & 0 & 0 && 0 & 0 & 1 & 0 && 0 & 0 & 0 & 0\\
0 & 0 & 0 & 0 && 0 & 0 & 0 & 0 && 0 & 0 & 0 & 1 && 0 & 0 & 0 & 0\\
\\
0 & 0 & 0 & 0 && 0 & 0 & 0 & 0 && 0 & 0 & 0 & 0 && 1 & 1 & 1 & 1\\
0 & 0 & 0 & 0 && 0 & 0 & 0 & 0 && 0 & 0 & 0 & 0 && 0 & 1 & 0 & 0\\
0 & 0 & 0 & 0 && 0 & 0 & 0 & 0 && 0 & 0 & 0 & 0 && 0 & 0 & 1 & 0\\
0 & 0 & 0 & 0 && 0 & 0 & 0 & 0 && 0 & 0 & 0 & 0 && 0 & 0 & 0 & 1\\
\end{array}
\right.
\left|
\begin{array}{r}
-20\\
10\\
10\\
10\\
\\
10\\
2\\
1\\
4\\
\\
10\\
3\\
2\\
1\\
\\
10\\
1\\
4\\
3\\
\end{array}
\right) .
\end{equation}
The solution vector x can be determined by inserting backwards.
\begin{equation*}
x  =  ( 1, 4, 3, 2, \  3, 2, 1, 4, \  4, 3, 2, 1, \  2, 1,  4, 3)
\end{equation*}
This results in the matrix corresponding to the Latin square in \autoref{lq4full} as the solution:
\begin{eqnarray*}
A & = &
\begin{pmatrix}
1 & 4 & 3 & 2\\
3 & 2 & 1 & 4\\
4 & 3 & 2 & 1\\
2 & 1 & 4 & 3\\
\end{pmatrix}.
\end{eqnarray*}

\subsection{Sudoku}
\subsubsection{Equations}
In a Sudoku with $a_{i,j} \in \{x \in \mathbb{N} \mid 1\leq x \leq n\}$ each row, each column and
each block has the sum:
\begin{equation*}
s := \sum_{k=1}^{n} k = \frac{n(n+1)}{2}.
\end{equation*}

This means that the following linear system of equations with $3n$ equations applies to a Sudoku of
order~$n= l \cdot m$ with blocks of size $l$ rows and $m$ columns:
\begin{eqnarray}
\sum_{i=1}^{n} a_{i,j} = s, &  j = 1, \ldots , n, & \mbox{equations $1, \ldots, n$ for $n$ rows}, \nonumber \\
\sum_{j=1}^{n} a_{i,j} = s, & i = 1, \ldots , n, & \mbox{equations $n + 1, \ldots, 2n$ for $n$ columns},  \label{e:LgsSdk} \\
\sum_{i=1}^{l}\left( \sum_{j=1}^{m} a_{p, q} \right) = s,  & k = 0, \ldots , (n-1),  &
     \mbox{equations $n + 2, \ldots, 3n$ for $n = l \cdot m $ blocks}, \nonumber \\
& & p = k - (k\bmod{l})+i, \nonumber \\
& & q = (k\bmod{l}) \cdot m+j. \nonumber
\end{eqnarray}
This system of equations corresponds to the linear equations of equation (2) of \cite{Arnold_Sudoku_2010}.
Overall, Sudokus in~\cite{Arnold_Sudoku_2010}
are described by a system of equations consisting of 1 linear and 2 nonlinear sub-systems of equations that
must be solved simultaneously, with all 3 sub-systems of equations having a common solution. In other
words, each solution of one of the 3 sub-systems of equations is also a solution of the other two
subsystems of equations. Therefore, a solution of the system of equations
(\ref{e:LgsSdk})) is also a solution of the other two
partial systems of equations. In \cite{Arnold_Sudoku_2010}, a Gröbner
basis is only given for Sudokus of order 4. Sudokus of a
higher order could not be solved for performance reasons.
The system of equations (\ref{e:LgsSdk})) describes the relationship
that each of the numbers 1 to n may occur at most
once per row, column or block. It does not contain any equations
that describe the relationship that each of
the numbers 1 to n must occur at least once per row and column.
This condition is irrelevant for the further
argumentation in this article.

\subsubsection{Transformation to Row reduced Echelon Form}\label{sec:lqzf}

In matrix notation, the system of equations is written as \( Bx = c \), where \( B = (b_{i,j}) \) is the coefficient matrix and \( c = (c_i) \) is the right-hand side vector.

The first \( 2n \) equations of a Sudoku system correspond to those of the Latin square (cf. equations~\ref{eq:lq:a}). This makes it clear that a Sudoku is essentially a Latin square with additional constraints.

For a Sudoku of order \( n = 4 \) with \( a_{i,j} \in \{1, 2, 3, 4\} \) and \( s = 10 \), the resulting linear system consists of the following 12 equations:

\begin{eqnarray*}
a_{1,1} + a_{2,1} + a_{3,1} + a_{4,1} & = & 10\\
a_{1,2} + a_{2,2} + a_{3,2} + a_{4,2} & = & 10\\
a_{1,3} + a_{2,3} + a_{3,3} + a_{4,3} & = & 10\\
a_{1,4} + a_{2,4} + a_{3,4} + a_{4,4} & = & 10\\[1ex]
a_{1,1} + a_{1,2} + a_{1,3} + a_{1,4} & = & 10\\
a_{2,1} + a_{2,2} + a_{2,3} + a_{2,4} & = & 10\\
a_{3,1} + a_{3,2} + a_{3,3} + a_{3,4} & = & 10\\
a_{4,1} + a_{4,2} + a_{4,3} + a_{4,4} & = & 10\\[1ex]
a_{1,1} + a_{1,2} + a_{2,1} + a_{2,2} & = & 10\\
a_{1,3} + a_{1,4} + a_{2,3} + a_{2,4} & = & 10\\
a_{3,1} + a_{3,2} + a_{4,1} + a_{4,2} & = & 10\\
a_{3,3} + a_{3,4} + a_{4,3} + a_{4,4} & = & 10
\end{eqnarray*}

In vector notation using variables \( x_i \), the general system can be written as:

\begin{eqnarray}
\sum_{i=0}^{n-1} x_{n \cdot i + k} = s, & \quad k = 1 , \ldots , n, & \quad \text{(Columns,)} \nonumber \\
\sum_{j=1}^{n} x_{n \cdot k + j} = s, & \quad k = 0 , \ldots , n-1, & \quad \text{(Rows,)} \label{eq:xSdk}\\
\sum^{l-1}_{j=0} \left( \sum^m_{i=1} x_{i+j \cdot n + p + q} \right) =  s, & \quad k = 0 , \ldots , n-1, & \quad \text{(Blocks)}\nonumber \\
& p = m \cdot k \bmod{l}, & \nonumber \\
& q = n \cdot (k - k \bmod{l}). & \nonumber
\end{eqnarray}

These equations correspond to Equation (2) in~\cite{Arnold_Sudoku_2010}. However, as already discussed in Equation~(\ref{e:LgsSdk}), they do not guarantee that each number from 1 to \( n \) appears \emph{exactly once} in every row and column.

For \( n = 4 \), the extended coefficient matrix is as follows:
\begin{equation*}
(B | c ) = \left(
\begin{array}{rrrrcrrrrcrrrrcrrrr}
		1 & 1 & 1 & 1 && 0 & 0 & 0 & 0 && 0 & 0 & 0 & 0 && 0 & 0 & 0 & 0\\
		0 & 0 & 0 & 0 && 1 & 1 & 1 & 1 && 0 & 0 & 0 & 0 && 0 & 0 & 0 & 0\\
		0 & 0 & 0 & 0 && 0 & 0 & 0 & 0 && 1 & 1 & 1 & 1 && 0 & 0 & 0 & 0\\
		0 & 0 & 0 & 0 && 0 & 0 & 0 & 0 && 0 & 0 & 0 & 0 && 1 & 1 & 1 & 1\\
		&\\
		1 & 0 & 0 & 0 && 1 & 0 & 0 & 0 && 1 & 0 & 0 & 0 && 1 & 0 & 0 & 0\\
		0 & 1 & 0 & 0 && 0 & 1 & 0 & 0 && 0 & 1 & 0 & 0 && 0 & 1 & 0 & 0\\
		0 & 0 & 1 & 0 && 0 & 0 & 1 & 0 && 0 & 0 & 1 & 0 && 0 & 0 & 1 & 0\\
		0 & 0 & 0 & 1 && 0 & 0 & 0 & 1 && 0 & 0 & 0 & 1 && 0 & 0 & 0 & 1\\
		&\\
		1 & 1 & 0 & 0 && 1 & 1 & 0 & 0 && 0 & 0 & 0 & 0 && 0 & 0 & 0 & 0\\
		0 & 0 & 1 & 1 && 0 & 0 & 1 & 1 && 0 & 0 & 0 & 0 && 0 & 0 & 0 & 0\\
		0 & 0 & 0 & 0 && 0 & 0 & 0 & 0 && 1 & 1 & 0 & 0 && 1 & 1 & 0 & 0\\
		0 & 0 & 0 & 0 && 0 & 0 & 0 & 0 && 0 & 0 & 1 & 1 && 0 & 0 & 1 & 1\\
\end{array}
\right.
\left|
\begin{array}{r}
10\\
10\\
10\\
10\\
\\
10\\
10\\
10\\
10\\
\\
10\\
10\\
10\\
10\\
\end{array}
\right).
\end{equation*}

This matrix was reduced to \emph{reduced row echelon form} using GNU Octave~\cite{gnu:octave:8.4}. The source code and results are available on GitHub~\cite{rp:gh:lq_sdk:lgs}.

The result is the transformed matrix:
\begin{equation}
(B | c ) = \left(
\begin{array}{rrrrcrrrrcrrrrcrrrr}
    1 & 0 & 0 & 0 && 0 & -1 & -1 & -1 && 0 & -1 & 0 & 0 && 0 & -1 & 0 & 0\\
    0 & 1 & 0 & 0 && 0 &  1 & 0 & 0 && 0 & 1 & 0 & 0 && 0 & 1 & 0 & 0\\
    0 & 0 & 1 & 0 && 0 & 0 & 1 & 0 && 0 & 0 & 0 & -1 && 0 & 0 & 0 & -1\\
    0 & 0 & 0 & 1 && 0 & 0 & 0 & 1 && 0 & 0 & 0 & 1 && 0 & 0 & 0 & 1\\
    &\\
    0 & 0 & 0 & 0 && 1 & 1 & 1 & 1 && 0 & 0 & 0 & 0 && 0 & 0 & 0 & 0\\
    0 & 0 & 0 & 0 && 0 & 0 & 0 & 0 && 1 & 1 & 0 & 0 && 0 & 0 & -1 & -1\\
    0 & 0 & 0 & 0 && 0 & 0 & 0 & 0 && 0 & 0 & 1 & 1 && 0 & 0 & 1 & 1\\
    0 & 0 & 0 & 0 && 0 & 0 & 0 & 0 && 0 & 0 & 0 & 0 && 1 & 1 & 1 & 1\\
    &\\
    0 & 0 & 0 & 0 && 0 & 0 & 0 & 0 && 0 & 0 & 0 & 0 && 0 & 0 & 0 & 0\\
    0 & 0 & 0 & 0 && 0 & 0 & 0 & 0 && 0 & 0 & 0 & 0 && 0 & 0 & 0 & 0\\
    0 & 0 & 0 & 0 && 0 & 0 & 0 & 0 && 0 & 0 & 0 & 0 && 0 & 0 & 0 & 0\\
    0 & 0 & 0 & 0 && 0 & 0 & 0 & 0 && 0 & 0 & 0 & 0 && 0 & 0 & 0 & 0\\
\end{array}
\right.
\left|
\begin{array}{r}
-10\\
10\\
0\\
10\\
\\
10\\
0\\
10\\
10\\
\\
0\\
0\\
0\\
0\\
\end{array}
\right). \label{sdk4koeff}
\end{equation}
The rank of the matrix can now be directly read off: \( \operatorname{rank}(B) = 8 \).

Applying the full Sudoku solution from Figure~\ref{sdk4full}, removing the pivot elements from the non-zero rows, and interpreting
them as unknowns yields the puzzle shown in Figure~\ref{sdk4rang}.
\begin{figure}[h]
\begin{tikzpicture}[scale=0.8]
\node at (-5, 0) { };
\draw[help lines] (0,0) grid (4,4);
\draw[help lines] (6,0) grid (10,4);
    \draw[line width=1.5pt]  (0,0) -- (4,0) -- (4,4) -- (0,4) -- (0,0) -- (4,0);
\draw[line width=1.5pt]  (0,2) -- (4,2);
\draw[line width=1.5pt]  (2,0) -- (2,4);
\node  at (0.5, 3.5) {$\not 1$};
\node  at (2.5, 2.5) {1};
\node  at (3.5, 1.5) {1};
\node  at (1.5, 0.5) {1};
\node  at (3.5, 3.5) {$\not 2$};
\node  at (1.5, 2.5) {2};
\node  at (2.5, 1.5) {$\not 2$};
\node  at (0.5, 0.5) {$\not 2$};
\node  at (2.5, 3.5) {$\not 3$};
\node  at (0.5, 2.5) {$\not 3$};
\node  at (1.5, 1.5) {3};
\node  at (3.5, 0.5) {3};
\node  at (1.5, 3.5) {$\not 4$};
\node  at (3.5, 2.5) {4};
\node  at (0.5, 1.5) {$\not 4$};
\node  at (2.5, 0.5) {4};
\draw[help lines] (6,0) grid (10,4);
    \draw[line width=1.5pt]  (6,0) -- (10,0) -- (10,4) -- (6,4) -- (6,0) -- (10,0);
\draw[line width=1.5pt]  (6,2) -- (10,2);
\draw[line width=1.5pt]  (8,0) -- (8,4);
\node  at (8.5, 2.5) {1};
\node  at (9.5, 1.5) {1};
\node  at (7.5, 0.5) {1};
\node  at (7.5, 2.5) {2};
\node  at (7.5, 1.5) {3};
\node  at (9.5, 0.5) {3};
\node  at (9.5, 2.5) {4};
\node  at (8.5, 0.5) {4};
\end{tikzpicture}
\caption{Example creation of a Sudoku of order $n = 4, l = 2, m = 2$  \label{sdk4rang}}
\end{figure}

If additional equations for the known values are added to the system~(\ref{sdk4koeff}) and the matrix is again brought into row echelon form, the resulting matrix has \emph{full rank}.

Using back-substitution, the solution vector \( x \) can then be uniquely determined.
\begin{equation*}
x  = ( 1, 4, 3, 2, \ 3, 2, 1, 4, \  4, 3, 2, 1, \  2, 1,  4, 3)
\end{equation*}

This yields the matrix corresponding to the Sudoku in \autoref{sdk4full} as the solution:
\begin{equation*}
A  =
\begin{pmatrix}
1 & 4 & 3 & 2\\
3 & 2 & 1 & 4\\
4 & 3 & 2 & 1\\
2 & 1 & 4 & 3\\
\end{pmatrix}.
\end{equation*}

\subsubsection{Examples for the Rank of the Coefficient Matrix for Several Sudokus}
Using GNU Octave, the coefficient matrix for Sudokus of orders~$4, 6, 8, 9$ was computed, and the row echelon form
and the rank~$rk$ were determined. The program used and the results are available on GitHub~\cite{rp:gh:lq_sdk:lgs}.

Based on the formula for the rank of Latin squares, the following formula for the rank $\overline{rk}$ is assumed:
\begin{equation}
\overline{rk} = 2n -1 + (l -1) * (m - 1) \label{rksdk}.
\end{equation}

The results for the rank $rk$ and $\overline{rk}$ are listed in the following table.
\begin{equation*}
\begin{tabular}{|r|r|r|r|r|}
\hline
$n$ & $l$ & $m$ & $rk$ & $\overline{rk}$\\
\hline
4 & 2 & 2 & 8 & 8\\
\hline
6 & 2 & 3 & 13 & 13\\
\hline
6 & 3 & 2 & 13 & 13\\
\hline
8 & 2 & 4 & 18 & 18\\
\hline
8 & 4 & 2 & 18 & 18\\
\hline
9 & 3 & 3 & 21 & 21\\
\hline
\end{tabular}
\end{equation*}

This listing provides only numerical evidence for the assumed formula for the rank. To derive
formula~(\ref{rksdk}), Sudokus are considered in which the pivot columns of the
coefficient matrix in row reduced echelon form are chosen as unknowns.
These are illustrated in a simplified graphical representation.
Empty fields represent unknowns, and known variables are denoted by `*'.
We show the distribution for Sudokus of orders: 4, 6, 8, 9
in figures \ref{unknownssdk4} to \ref{unknownssdk9}.

\begin{figure}
\begin{center}
\begin{tikzpicture}[scale=0.5]
\draw[help lines] (0,0) grid (4,4);
\draw[line width=1.5pt]  (0,0) -- (4,0) -- (4,4) --(0,4) -- (0,0) -- (4,0);
\draw[line width=1.5pt]  (0,2) -- (4,2);
\draw[line width=1.5pt]  (2,0) -- (2,4);
\node  at (2.5, 2.5) {*};
\node  at (3.5, 1.5) {*};
\node  at (1.5, 0.5) {*};
\node  at (1.5, 2.5) {*};
\node  at (1.5, 1.5) {*};
\node  at (3.5, 0.5) {*};
\node  at (3.5, 2.5) {*};
\node  at (2.5, 0.5) {*};
\end{tikzpicture}
\end{center}
\caption{Distribution of unknowns for a Sudoku of order 4} \label{unknownssdk4}
\end{figure}
\begin{figure}
\begin{center}
\begin{tikzpicture}[scale=0.5]
\draw[help lines] (0,0) grid (6,6);
\draw[line width=1.5pt]  (0,0) -- (6,0) -- (6,6) -- (0,6) -- (0,0) -- (6,0);
\draw[line width=1.5pt]  (0,2) -- (6,2);
\draw[line width=1.5pt]  (0,4) -- (6,4);
\draw[line width=1.5pt]  (3,0) -- (3,6);
\node  at (1.5, 4.5) {*};
\node  at (2.5, 4.5) {*};
\node  at (3.5, 4.5) {*};
\node  at (4.5, 4.5) {*};
\node  at (5.5, 4.5) {*};
\node  at (1.5, 3.5) {*};
\node  at (2.5, 3.5) {*};
\node  at (4.5, 3.5) {*};
\node  at (5.5, 3.5) {*};
\node  at (1.5, 2.5) {*};
\node  at (2.5, 2.5) {*};
\node  at (3.5, 2.5) {*};
\node  at (4.5, 2.5) {*};
\node  at (5.5, 2.5) {*};
\node  at (1.5, 1.5) {*};
\node  at (2.5, 1.5) {*};
\node  at (4.5, 1.5) {*};
\node  at (5.5, 1.5) {*};
\node  at (1.5, 0.5) {*};
\node  at (2.5, 0.5) {*};
\node  at (3.5, 0.5) {*};
\node  at (4.5, 0.5) {*};
\node  at (5.5, 0.5) {*};
\draw[help lines] (8,0) grid (14,6);
\draw[line width=1.5pt]  (8,0) -- (14,0) -- (14,6) -- (8,6) -- (8,0) -- (14, 0);
\draw[line width=1.5pt]  (8,3) -- (14,3);
\draw[line width=1.5pt]  (10,0) -- (10,6);
\draw[line width=1.5pt]  (12,0) -- (12,6);
\node  at (9.5, 4.5) {*};
\node  at (10.5, 4.5) {*};
\node  at (11.5, 4.5) {*};
\node  at (12.5, 4.5) {*};
\node  at (13.5, 4.5) {*};
\node  at (9.5, 3.5) {*};
\node  at (10.5, 3.5) {*};
\node  at (11.5, 3.5) {*};
\node  at (12.5, 3.5) {*};
\node  at (13.5, 3.5) {*};
\node  at (9.5, 2.5) {*};
\node  at (11.5, 2.5) {*};
\node  at (13.5, 2.5) {*};
\node  at (9.5, 1.5) {*};
\node  at (10.5, 1.5) {*};
\node  at (11.5, 1.5) {*};
\node  at (12.5, 1.5) {*};
\node  at (13.5, 1.5) {*};
\node  at (9.5, 0.5) {*};
\node  at (10.5, 0.5) {*};
\node  at (11.5, 0.5) {*};
\node  at (12.5, 0.5) {*};
\node  at (13.5, 0.5) {*};
\end{tikzpicture}
\end{center}
\caption{Distribution of unknowns for a Sudoku of order 6}
\end{figure}

\begin{figure}
\begin{center}
\begin{tikzpicture}[scale=0.5]
\draw[help lines] (0,0) grid (8,8);
\draw[line width=1.5pt]  (0,0) -- (8,0) -- (8,8) --(0,8) -- (0,0) -- (8,0);
\draw[line width=1.5pt]  (0,2) -- (8,2);
\draw[line width=1.5pt]  (0,4) -- (8,4);
\draw[line width=1.5pt]  (0,6) -- (8,6);
\draw[line width=1.5pt]  (4,0) -- (4,8);
\node  at (1.5, 6.5) {*};
\node  at (2.5, 6.5) {*};
\node  at (3.5, 6.5) {*};
\node  at (4.5, 6.5) {*};
\node  at (5.5, 6.5) {*};
\node  at (6.5, 6.5) {*};
\node  at (7.5, 6.5) {*};
\node  at (1.5, 5.5) {*};
\node  at (2.5, 5.5) {*};
\node  at (3.5, 5.5) {*};
\node  at (5.5, 5.5) {*};
\node  at (6.5, 5.5) {*};
\node  at (7.5, 5.5) {*};
\node  at (1.5, 4.5) {*};
\node  at (2.5, 4.5) {*};
\node  at (3.5, 4.5) {*};
\node  at (4.5, 4.5) {*};
\node  at (5.5, 4.5) {*};
\node  at (6.5, 4.5) {*};
\node  at (7.5, 4.5) {*};
\node  at (1.5, 3.5) {*};
\node  at (2.5, 3.5) {*};
\node  at (3.5, 3.5) {*};
\node  at (5.5, 3.5) {*};
\node  at (6.5, 3.5) {*};
\node  at (7.5, 3.5) {*};
\node  at (1.5, 2.5) {*};
\node  at (2.5, 2.5) {*};
\node  at (3.5, 2.5) {*};
\node  at (4.5, 2.5) {*};
\node  at (5.5, 2.5) {*};
\node  at (6.5, 2.5) {*};
\node  at (7.5, 2.5) {*};
\node  at (1.5, 1.5) {*};
\node  at (2.5, 1.5) {*};
\node  at (3.5, 1.5) {*};
\node  at (5.5, 1.5) {*};
\node  at (6.5, 1.5) {*};
\node  at (7.5, 1.5) {*};
\node  at (1.5, 0.5) {*};
\node  at (2.5, 0.5) {*};
\node  at (3.5, 0.5) {*};
\node  at (4.5, 0.5) {*};
\node  at (5.5, 0.5) {*};
\node  at (6.5, 0.5) {*};
\node  at (7.5, 0.5) {*};
\draw[help lines] (10,0) grid (18,8);
\draw[line width=1.5pt]  (10,0) -- (18,0) -- (18,8) -- (10,8) -- (10,0) -- (18,0);
\draw[line width=1.5pt]  (10,4) -- (18,4);
\draw[line width=1.5pt]  (12,0) -- (12,8);
\draw[line width=1.5pt]  (14,0) -- (14,8);
\draw[line width=1.5pt]  (16,0) -- (16,8);
\node  at (11.5, 6.5) {*};
\node  at (12.5, 6.5) {*};
\node  at (13.5, 6.5) {*};
\node  at (14.5, 6.5) {*};
\node  at (15.5, 6.5) {*};
\node  at (16.5, 6.5) {*};
\node  at (17.5, 6.5) {*};
\node  at (11.5, 5.5) {*};
\node  at (12.5, 5.5) {*};
\node  at (13.5, 5.5) {*};
\node  at (14.5, 5.5) {*};
\node  at (15.5, 5.5) {*};
\node  at (16.5, 5.5) {*};
\node  at (17.5, 5.5) {*};
\node  at (11.5, 4.5) {*};
\node  at (12.5, 4.5) {*};
\node  at (13.5, 4.5) {*};
\node  at (14.5, 4.5) {*};
\node  at (15.5, 4.5) {*};
\node  at (16.5, 4.5) {*};
\node  at (17.5, 4.5) {*};
\node  at (11.5, 3.5) {*};
\node  at (13.5, 3.5) {*};
\node  at (15.5, 3.5) {*};
\node  at (17.5, 3.5) {*};
\node  at (11.5, 2.5) {*};
\node  at (12.5, 2.5) {*};
\node  at (13.5, 2.5) {*};
\node  at (14.5, 2.5) {*};
\node  at (15.5, 2.5) {*};
\node  at (16.5, 2.5) {*};
\node  at (17.5, 2.5) {*};
\node  at (11.5, 1.5) {*};
\node  at (12.5, 1.5) {*};
\node  at (13.5, 1.5) {*};
\node  at (14.5, 1.5) {*};
\node  at (15.5, 1.5) {*};
\node  at (16.5, 1.5) {*};
\node  at (17.5, 1.5) {*};
\node  at (11.5, 0.5) {*};
\node  at (12.5, 0.5) {*};
\node  at (13.5, 0.5) {*};
\node  at (14.5, 0.5) {*};
\node  at (15.5, 0.5) {*};
\node  at (16.5, 0.5) {*};
\node  at (17.5, 0.5) {*};
\end{tikzpicture}
\end{center}
\caption{Distribution of unknowns for a Sudoku of order 8}
\end{figure}

\begin{figure}
\begin{center}
\begin{tikzpicture}[scale=0.5]
\node at (-2, 0) { };
\draw[help lines] (0,0) grid (9,9);
    \draw[line width=1.5pt]  (0,0) -- (9,0) -- (9,9) --(0,9) -- (0,0) -- (9,0);
\draw[line width=1.5pt]  (0,3) -- (9,3);
\draw[line width=1.5pt]  (0,6) -- (9,6);
\draw[line width=1.5pt]  (3,0) -- (3,9);
\draw[line width=1.5pt]  (6,0) -- (6,9);
\node  at (1.5, 7.5) {*};
\node  at (2.5, 7.5) {*};
\node  at (3.5, 7.5) {*};
\node  at (4.5, 7.5) {*};
\node  at (5.5, 7.5) {*};
\node  at (6.5, 7.5) {*};
\node  at (7.5, 7.5) {*};
\node  at (8.5, 7.5) {*};
\node  at (1.5, 6.5) {*};
\node  at (2.5, 6.5) {*};
\node  at (3.5, 6.5) {*};
\node  at (4.5, 6.5) {*};
\node  at (5.5, 6.5) {*};
\node  at (6.5, 6.5) {*};
\node  at (7.5, 6.5) {*};
\node  at (8.5, 6.5) {*};
\node  at (1.5, 5.5) {*};
\node  at (2.5, 5.5) {*};
\node  at (4.5, 5.5) {*};
\node  at (5.5, 5.5) {*};
\node  at (7.5, 5.5) {*};
\node  at (8.5, 5.5) {*};
\node  at (1.5, 4.5) {*};
\node  at (2.5, 4.5) {*};
\node  at (3.5, 4.5) {*};
\node  at (4.5, 4.5) {*};
\node  at (5.5, 4.5) {*};
\node  at (6.5, 4.5) {*};
\node  at (7.5, 4.5) {*};
\node  at (8.5, 4.5) {*};
\node  at (1.5, 3.5) {*};
\node  at (2.5, 3.5) {*};
\node  at (3.5, 3.5) {*};
\node  at (4.5, 3.5) {*};
\node  at (5.5, 3.5) {*};
\node  at (6.5, 3.5) {*};
\node  at (7.5, 3.5) {*};
\node  at (8.5, 3.5) {*};
\node  at (1.5, 2.5) {*};
\node  at (2.5, 2.5) {*};
\node  at (4.5, 2.5) {*};
\node  at (5.5, 2.5) {*};
\node  at (7.5, 2.5) {*};
\node  at (8.5, 2.5) {*};
\node  at (1.5, 1.5) {*};
\node  at (2.5, 1.5) {*};
\node  at (3.5, 1.5) {*};
\node  at (4.5, 1.5) {*};
\node  at (5.5, 1.5) {*};
\node  at (6.5, 1.5) {*};
\node  at (7.5, 1.5) {*};
\node  at (8.5, 1.5) {*};
\node  at (1.5, 0.5) {*};
\node  at (2.5, 0.5) {*};
\node  at (3.5, 0.5) {*};
\node  at (4.5, 0.5) {*};
\node  at (5.5, 0.5) {*};
\node  at (6.5, 0.5) {*};
\node  at (7.5, 0.5) {*};
\node  at (8.5, 0.5) {*};
\end{tikzpicture}
\end{center}
\caption{Distribution of unknowns for a Sudoku of order 9} \label{unknownssdk9}
\end{figure}

\subsubsection{Rank of the Coefficient Matrix for Sudoku}

In the graphical representation, a recurring pattern can be observed. As in Latin squares, the unknowns appear in the first row and the first column. Additionally, the elements in the first row and first column of each block that does not contain elements from the first row or first column also appear as unknowns.

This can be explained by the additional constraint in Sudoku regarding the blocks. For each complete block that does not contain any elements from the first row or the first column, an additional unknown can be added. This explains why the rank of the coefficient matrix in Sudoku is higher than that in Latin squares. Each Sudoku has $(l - 1) * (m - 1)$ blocks without elements from the first row and first column, and the rank increases exactly by this number.

For a Sudoku of order 9, this results in a rank of 21 for the linear system of equations. This exact value is also given in \cite{Demoen2014-bc} as the minimum number of logical constraints required for a Sudoku. Is this numerical equality a coincidence or mathematics – a deeper connection just waiting to be discovered?

\subsubsection{Rank of the Coefficient Matrix for Arbitrary Sudokus}
Based on the observations of the rank in the analyzed Sudokus, it can be shown that the rank of the coefficient matrix is given by:
\begin{equation}
rk = 2n - 1 + (l - 1) * (m - 1).
\end{equation}

If one considers the $3n$ equations from (\ref{eq:xSdk}) separately, it becomes apparent that each group of $n$ equations — for columns, rows, and blocks — is already in row echelon form.

For the equations corresponding to the columns
\begin{equation}
\sum_{i=0}^{n-1} x_{n \cdot i + k} = s, \quad k = 1 , \ldots , n
\end{equation}
the pivot elements have the indices:
\begin{equation*}
i = 1, \ldots, n.
\end{equation*}
If this is represented graphically, with the pivot elements as $x_i$ and other elements as $*$, we obtain the depiction shown in Figure~\ref{fig:Zfcols}.
\begin{figure}
\begin{center}
\begin{tikzpicture}[scale=0.7]
\node at (-2, 0) { };
\draw[help lines] (0,0) grid (4,4);
\draw[line width=1.5pt]  (0,0) -- (4,0) -- (4,4) -- (0,4) -- (0,0) -- (4,0);
\node  at (0.5, 3.5) {$x_1$};
\node  at (1.5, 3.5) {$x_2$};
\node  at (2.5, 3.5) {$\cdots$};
\node  at (3.5, 3.5) {$x_n$};
\node  at (0.5, 2.5) {*};
\node  at (1.5, 2.5) {*};
\node  at (2.5, 2.5) {$\cdots$};
\node  at (3.5, 2.5) {*};
\node  at (0.5, 1.6) {$\vdots$};
\node  at (1.5, 1.6) {$\vdots$};
\node  at (3.5, 1.6) {$\vdots$};
\node  at (0.5, 0.5) {*};
\node  at (1.5, 0.5) {*};
\node  at (2.5, 0.5) {$\cdots$};
\node  at (3.5, 0.5) {*};
\end{tikzpicture}
\end{center}
\caption{Pivot elements column equations}  \label{fig:Zfcols}
\end{figure}

For the equations corresponding to the rows
\begin{equation}
\sum_{j=1}^{n} x_{n \cdot k + j} = s, \quad k = 0 , \ldots , n-1
\end{equation}
the pivot elements have the indices:
\begin{equation*}
i = 1, n+1, \ldots, ( n - 1) * n + 1.
\end{equation*}
When represented graphically, using $x_{i_j}$ for the pivot elements and $*$ for the others, the resulting figure is shown in Figure~\ref{fig:Zfrows}.
\begin{figure}
\begin{center}
\begin{tikzpicture}[scale=0.7]
\node at (-1, 0) { };
\draw[help lines] (0,0) grid (4,4);
    \draw[line width=1.5pt]  (0,0) -- (4,0) -- (4,4) -- (0,4) -- (0,0) -- (4,0);
\node  at (0.5, 3.5) {$x_{i_{1}}$};
\node  at (1.5, 3.5) {*};
\node  at (2.5, 3.5) {$\cdots$};
\node  at (3.5, 3.5) {*};
\node  at (0.5, 2.5) {$x_{i_{2}}$};
\node  at (1.5, 2.5) {*};
\node  at (2.5, 2.5) {$\cdots$};
\node  at (3.5, 2.5) {*};
\node  at (0.5, 1.6) {$\vdots$};
\node  at (1.5, 1.6) {$\vdots$};
\node  at (3.5, 1.6) {$\vdots$};
\node  at (0.5, 0.5) {$x_{i_{n}}$};
\node  at (1.5, 0.5) {*};
\node  at (2.5, 0.5) {$\cdots$};
\node  at (3.5, 0.5) {*};
\node  at (9.0, 0.5) {with $i_k = 1, n+1, \ldots, (n-1) * n + 1$};
\end{tikzpicture}
\end{center}
\caption{Pivot elements row equations}  \label{fig:Zfrows}
\end{figure}

For the equations corresponding to blocks of size $l$ rows and $m$ columns (cf. \autoref{eq:xSdk})
\begin{eqnarray}
\sum^{l-1}_{j=0} \left( \sum^m_{i=1} x_{i+j*n+p+q} \right) =  s & &  k = 0 , \ldots , n-1,  \label{eq:k} \\
& & p = m \cdot k \bmod{l}, \nonumber \\
& & q = n \cdot (k - k \bmod{l}). \nonumber
\end{eqnarray}
the pivot elements $x_{r_h}$ have the indices:
\begin{equation}
\begin{array}{rcrrcl}
    r_h & = & 1, & m+1, & \ldots, & (n-l) * n + ( l - 1) * m + 1, \\
    h & = & 1, & 2, & \ldots, & n
\end{array} \label{eq:r}
\end{equation}
since $l \cdot m = n$.
When represented graphically with the pivot elements labeled as $x_{r_h}$ and the others as~$*$,
the depiction is shown in Figure~\ref{fig:Zfblocks}.
\begin{figure}
\begin{center}
\begin{tikzpicture}[scale=0.7]
\draw[help lines] (0,0) grid (16,16);
\draw[line width=1.5pt]  (0,0) -- (16,0) -- (16,16) -- (0,16) -- (0,0) -- (16,0);
\draw[line width=1.5pt]  (4,0) -- (4,16);
\draw[line width=1.5pt]  (8,0) -- (8,16);
\draw[line width=1.5pt]  (12,0) -- (12,16);
\draw[line width=1.5pt]  (0,4) -- (16,4);
\draw[line width=1.5pt]  (0,8) -- (16,8);
\draw[line width=1.5pt]  (0,12) -- (16,12);
%  row 1 of blocks
\node  at (0.5, 15.5) {$x_{r_1}$};
\node  at (1.5, 15.5) {*};
\node  at (2.5, 15.5) {$\cdots$};
\node  at (3.5, 15.5) {*};
\node  at (4.5, 15.5) {$x_{r_2}$};
\node  at (5.5, 15.5) {*};
\node  at (6.5, 15.5) {$\cdots$};
\node  at (7.5, 15.5) {*};
\node  at (8.5, 15.5) {$\cdots$};
\node  at (9.5, 15.5) {$\cdots$};
\node  at (10.5, 15.5) {$\cdots$};
\node  at (11.5, 15.5) {$\cdots$};
\node  at (12.5, 15.5) {$x_{r_{l}}$};
\node  at (13.5, 15.5) {*};
\node  at (14.5, 15.5) {$\cdots$};
\node  at (15.5, 15.5) {*};
\node  at (0.5, 14.5) {*};
\node  at (1.5, 14.5) {$\cdots$};
\node  at (2.5, 14.5) {$\cdots$};
\node  at (3.5, 14.5) {*};
\node  at (4.5, 14.5) {*};
\node  at (5.5, 14.5) {$\cdots$};
\node  at (6.5, 14.5) {$\cdots$};
\node  at (7.5, 14.5) {*};
\node  at (8.5, 14.5) {$\cdots$};
\node  at (9.5, 14.5) {$\cdots$};
\node  at (10.5, 14.5) {$\cdots$};
\node  at (11.5, 14.5) {$\cdots$};
\node  at (12.5, 14.5) {*};
\node  at (13.5, 14.5) {$\cdots$};
\node  at (14.5, 14.5) {$\cdots$};
\node  at (15.5, 14.5) {*};
\node  at (0.5, 13.5) {$\vdots$};
\node  at (3.5, 13.5) {$\vdots$};
\node  at (4.5, 13.5) {$\vdots$};
\node  at (7.5, 13.5) {$\vdots$};
\node  at (12.5, 13.5) {$\vdots$};
\node  at (15.5, 13.5) {$\vdots$};
\node  at (0.5, 12.5) {*};
\node  at (1.5, 12.5) {$\cdots$};
\node  at (2.5, 12.5) {$\cdots$};
\node  at (3.5, 12.5) {*};
\node  at (4.5, 12.5) {*};
\node  at (5.5, 12.5) {$\cdots$};
\node  at (6.5, 12.5) {$\cdots$};
\node  at (7.5, 12.5) {*};
\node  at (8.5, 12.5) {$\cdots$};
\node  at (9.5, 12.5) {$\cdots$};
\node  at (10.5, 12.5) {$\cdots$};
\node  at (11.5, 12.5) {$\cdots$};
\node  at (12.5, 12.5) {*};
\node  at (13.5, 12.5) {$\cdots$};
\node  at (14.5, 12.5) {$\cdots$};
\node  at (15.5, 12.5) {*};
% 2nd row of blocks
\node  at (0.5, 11.5) {$x_{r_{h}}$};
\node  at (1.5, 11.5) {*};
\node  at (2.5, 11.5) {$\cdots$};
\node  at (3.5, 11.5) {*};
\node  at (4.5, 11.5) {$x_{r_{h}}$};
\node  at (5.5, 11.5) {*};
\node  at (6.5, 11.5) {$\cdots$};
\node  at (7.5, 11.5) {*};
\node  at (8.5, 11.5) {$\cdots$};
\node  at (9.5, 11.5) {$\cdots$};
\node  at (10.5, 11.5) {$\cdots$};
\node  at (11.5, 11.5) {$\cdots$};
\node  at (12.5, 11.5) {$x_{r_{h}}$};
\node  at (13.5, 11.5) {*};
\node  at (14.5, 11.5) {$\cdots$};
\node  at (15.5, 11.5) {*};
\node  at (0.5, 10.5) {*};
\node  at (1.5, 10.5) {$\cdots$};
\node  at (2.5, 10.5) {$\cdots$};
\node  at (3.5, 10.5) {*};
\node  at (4.5, 10.5) {*};
\node  at (5.5, 10.5) {$\cdots$};
\node  at (6.5, 10.5) {$\cdots$};
\node  at (7.5, 10.5) {*};
\node  at (8.5, 10.5) {$\cdots$};
\node  at (9.5, 10.5) {$\cdots$};
\node  at (10.5, 10.5) {$\cdots$};
\node  at (11.5, 10.5) {$\cdots$};
\node  at (12.5, 10.5) {*};
\node  at (13.5, 10.5) {$\cdots$};
\node  at (14.5, 10.5) {$\cdots$};
\node  at (15.5, 10.5) {*};
\node  at (0.5, 9.5) {$\vdots$};
\node  at (3.5, 9.5) {$\vdots$};
\node  at (4.5, 9.5) {$\vdots$};
\node  at (7.5, 9.5) {$\vdots$};
\node  at (12.5, 9.5) {$\vdots$};
\node  at (15.5, 9.5) {$\vdots$};
\node  at (0.5, 8.5) {*};
\node  at (1.5, 8.5) {$\cdots$};
\node  at (2.5, 8.5) {$\cdots$};
\node  at (3.5, 8.5) {*};
\node  at (4.5, 8.5) {*};
\node  at (5.5, 8.5) {$\cdots$};
\node  at (6.5, 8.5) {$\cdots$};
\node  at (7.5, 8.5) {*};
\node  at (8.5, 8.5) {$\cdots$};
\node  at (9.5, 8.5) {$\cdots$};
\node  at (10.5, 8.5) {$\cdots$};
\node  at (11.5, 8.5) {$\cdots$};
\node  at (12.5, 8.5) {*};
\node  at (13.5, 8.5) {$\cdots$};
\node  at (14.5, 8.5) {$\cdots$};
\node  at (15.5, 8.5) {*};
% third row
\node  at (0.5, 7.5) {$\vdots$};
\node  at (3.5, 7.5) {$\vdots$};
\node  at (4.5, 7.5) {$\vdots$};
\node  at (7.5, 7.5) {$\vdots$};
\node  at (12.5, 7.5) {$\vdots$};
\node  at (15.5, 7.5) {$\vdots$};
\node  at (0.5, 6.5) {$\vdots$};
\node  at (3.5, 6.5) {$\vdots$};
\node  at (4.5, 6.5) {$\vdots$};
\node  at (7.5, 6.5) {$\vdots$};
\node  at (12.5, 6.5) {$\vdots$};
\node  at (15.5, 6.5) {$\vdots$};
\node  at (0.5, 5.5) {$\vdots$};
\node  at (3.5, 5.5) {$\vdots$};
\node  at (4.5, 5.5) {$\vdots$};
\node  at (7.5, 5.5) {$\vdots$};
\node  at (12.5, 5.5) {$\vdots$};
\node  at (15.5, 5.5) {$\vdots$};
\node  at (0.5, 4.5) {$\vdots$};
\node  at (3.5, 4.5) {$\vdots$};
\node  at (4.5, 4.5) {$\vdots$};
\node  at (7.5, 4.5) {$\vdots$};
\node  at (12.5, 4.5) {$\vdots$};
\node  at (15.5, 4.5) {$\vdots$};
% last row of blocks
\node  at (0.5, 3.5) {$x_{r_{h}}$};
\node  at (1.5, 3.5) {*};
\node  at (2.5, 3.5) {$\cdots$};
\node  at (3.5, 3.5) {*};
\node  at (4.5, 3.5) {$x_{r_{h}}$};
\node  at (5.5, 3.5) {*};
\node  at (6.5, 3.5) {$\cdots$};
\node  at (7.5, 3.5) {*};
\node  at (8.5, 3.5) {$\cdots$};
\node  at (9.5, 3.5) {$\cdots$};
\node  at (10.5, 3.5) {$\cdots$};
\node  at (11.5, 3.5) {$\cdots$};
\node  at (12.5, 3.5) {$x_{r_{n}}$};
\node  at (13.5, 3.5) {*};
\node  at (14.5, 3.5) {$\cdots$};
\node  at (15.5, 3.5) {*};
\node  at (0.5, 2.5) {*};
\node  at (1.5, 2.5) {$\cdots$};
\node  at (2.5, 2.5) {$\cdots$};
\node  at (3.5, 2.5) {*};
\node  at (4.5, 2.5) {*};
\node  at (5.5, 2.5) {$\cdots$};
\node  at (6.5, 2.5) {$\cdots$};
\node  at (7.5, 2.5) {*};
\node  at (8.5, 2.5) {$\cdots$};
\node  at (9.5, 2.5) {$\cdots$};
\node  at (10.5, 2.5) {$\cdots$};
\node  at (11.5, 2.5) {$\cdots$};
\node  at (12.5, 2.5) {*};
\node  at (13.5, 2.5) {$\cdots$};
\node  at (14.5, 2.5) {$\cdots$};
\node  at (15.5, 2.5) {*};
\node  at (0.5, 1.5) {$\vdots$};
\node  at (3.5, 1.5) {$\vdots$};
\node  at (4.5, 1.5) {$\vdots$};
\node  at (7.5, 1.5) {$\vdots$};
\node  at (12.5, 1.5) {$\vdots$};
\node  at (15.5, 1.5) {$\vdots$};
\node  at (0.5, 0.5) {*};
\node  at (1.5, 0.5) {$\cdots$};
\node  at (2.5, 0.5) {$\cdots$};
\node  at (3.5, 0.5) {*};
\node  at (4.5, 0.5) {*};
\node  at (5.5, 0.5) {$\cdots$};
\node  at (6.5, 0.5) {$\cdots$};
\node  at (7.5, 0.5) {*};
\node  at (8.5, 0.5) {$\cdots$};
\node  at (9.5, 0.5) {$\cdots$};
\node  at (10.5, 0.5) {$\cdots$};
\node  at (11.5, 0.5) {$\cdots$};
\node  at (12.5, 0.5) {*};
\node  at (13.5, 0.5) {$\cdots$};
\node  at (14.5, 0.5) {$\cdots$};
\node  at (15.5, 0.5) {*};
\end{tikzpicture}
\end{center}
\caption{Pivot elements block equations} \label{fig:Zfblocks}
\end{figure}

Comparing the representations of the pivot elements, one observes that some pivot elements appear in only one of the representations, while others are present in multiple representations. The pivot variables that appear multiple times are:

\begin{enumerate}
\item $x_1$ is a pivot variable in all three representations, and thus in all three systems of equations.
\item $x_i$ with $i = m+1, \ldots , (l-1) * m + 1$ are pivot variables in the representations of the equations for the columns and the equations for the blocks.
\item $x_i$ with $i = l + 1, \ldots, (m - 1) * l + 1$ are pivot variables in the representations of the equations for the rows and the equations for the blocks.
\end{enumerate}

Let $z = 1, \ldots, 3n$ be the index of the equations in the system (\ref{eq:xSdk}). Selecting one equation for each pivot variable yields a system of equations in row echelon form.

One possible selection, without loss of generality, consists of the $n$ equations for the columns, i.e.,
\begin{equation*}
z = 1, \ldots, n,
\end{equation*}
the $n - 1$ equations for rows $2, \ldots, n$, i.e.,
\begin{equation*}
z = n+1, \ldots, (n - 1) * n + 1,
\end{equation*}
and the $(l - 1) * (m - 1)$ equations for the blocks that do not have a pivot variable in row 1 or column 1, i.e.,
\begin{eqnarray*}
& & z = n + l + 1, \ldots, (m - l) * l + 1 \\
& \wedge & z \neq n+1, \ldots, (n - 1) * n + 1.
\end{eqnarray*}

From the number of equations in this system in row echelon form, we derive the minimal rank as:
\begin{equation*}
rk \leq 2n - 1 + (l - 1) * (m - 1).
\end{equation*}

This is also the maximum rank. If it were not the maximum rank, then there would exist a row in the linear system
(\ref{e:LgsSdk})
with a pivot element at a position marked with a `$*$' in Figure \ref{fig:Zfblocks}.

\qed

Analogous to Latin squares, the same applies to Sudoku: depending on the choice of 4 unknowns in 2 rows, 2 columns, and 2 blocks, there may be one or two solutions. An example can be found in \cite{Delahaye2006}.

\section{Conclusion}

This article demonstrates that Latin squares and Sudokus of arbitrary order can be described and uniquely solved using
a linear system of equations. This confirms the hypothesis that back-substitution can be applied when solving Sudokus.
Thus, we have proved the existence of partial Latin squares and Sudokus that can be solved in polynomial time.
The presented algorithm belongs to complexity class P.
Usually, completing Latin squares and Sudokus is considered a problem of
complexity class NP (nondeterministic polynomial time) \cite{Colbourn1984,Thom2007}.
It may be that Latin squares and Sudokus are examples of systems in which
Schaefer’s dichotomy of complexity is applicable.

The rank of this system of equations is a function of the order of the Latin square or Sudoku.
Interestingly, the rank determined here exactly matches the minimal number
of logical constraints necessary to solve a Sudoku using pure logic.

However, the rank of the system is significantly smaller than the number
of unknowns in most Sudokus published in magazines.
The reason is that the linear system only describes the necessary conditions for a Latin square or Sudoku.
Additional equations are required to express the full set of constraints.
Arnold~\cite{Arnold_Sudoku_2010} describes these equations as a sum–product system.
The linear system presented here corresponds only to the “sum” part of these equations,
while the “product” part leads to a nonlinear system of equations.
This provides a starting point for further investigation.

\newpage
\printbibliography

\end{document}